\documentclass[a4paper,11pt]{article}

\def\empha{\em}

\def\bemphas{}
\def\pemphas{}
\def\scsc{\sc}
\newtheorem{thm}{Theorem}[section]

\newtheorem{prop}[thm]{Proposition}
\newtheorem{defi}[thm]{Definition}

\long
\def\MSC#1\EndMSC{\def\arg{#1}\ifx\arg\empty\relax\else
      {\par\narrower\noindent%
      2000 Mathematics Subject Classification. #1\par}\fi}

\long
\def\KEY#1\EndKEY{\def\arg{#1}\ifx\arg\empty\relax\else
    {\par\narrower\noindent%
      Keywords and Phrases: #1\par}\fi}

\title
{Singular Poisson-K\"ahler geometry of certain adjoint quotients}

\author{Johannes Huebschmann \footnote{Support by the Deutsche
Forschungsgemeinschaft in the framework of a
Mercator-professorship is gratefully
acknowledged}\\
Universit\'e des
Sciences et Technologies de Lille\\ UFR de Math\'ematiques, CNRS-UMR 8524
\\ 59655 VILLENEUVE D'ASCQ, C\'edex, France
\\ and
\\ Institute
for Theoretical Physics\\ Universit\"at Leipzig\\
04109 LEIPZIG, Germany \\
{Johannes.Huebschmann@math.univ-lille1.fr} }

\usepackage{amssymb}
\usepackage[square,authoryear]{natbib}
\usepackage{amsmath}
\usepackage{amscd}

\numberwithin{equation}{section}

\begin{document}
\setcounter{page}{1}

\maketitle

\begin{abstract}
The K\"ahler quotient of a complex reductive Lie group relative to
the conjugation action carries a complex algebraic stratified
K\"ahler structure which reflects the geometry of the group. For
the group $\mathrm{SL}(n,\mathbb C)$, we interpret the resulting
singular Poisson-K\"ahler geometry of the quotient in terms of
complex discriminant varieties and variants thereof.
\end{abstract}

\MSC Primary:  17B63 17B65 17B66 17B81 53D17 53D20 53D50 70H45
81S10; Secondary: 14L24 14L30 32C20 32Q15 32S05 32S60 \EndMSC

\KEY Stratified symplectic space, complex analytic space, complex
algebraic space, real semialgebraic space, complex analytic
stratified K\"ahler space, complex algebraic stratified K\"ahler
space, K\"ahler reduction, constrained system, invariant theory,
adjoint quotient, discriminant variety, spherical pendulum
 \EndKEY

\section{Adjoint quotients}

Let $K$ be a compact Lie group,  let $K^{\mathbb C}$ be its
complexification, and let $\mathfrak k$ be the Lie algebra of $K$.
The polar map from $K \times \mathfrak k$ to $K^{\mathbb C}$  is
given by the assignment to $(x,Y) \in K \times \mathfrak k$ of
$x\, \mathrm{exp}(iY) \in K^{\mathbb C}$, and this map is well
known to be a $K$-biinvariant diffeomorphism. We endow the Lie
algebra $\mathfrak k$  with an invariant (positive definite) inner
product; by means of this inner product, we identify $\mathfrak k$
with its dual $\mathfrak k^*$ and, furthermore, the total space
$\mathrm TK$ of the tangent bundle of $K$  with the total space
$\mathrm T^*K$ of the cotangent bundle of $K$. The composite
\[
\mathrm TK \longrightarrow K \times \mathfrak k \longrightarrow
K^{\mathbb C}
\]
of the inverse $\mathrm TK \longrightarrow K \times \mathfrak k$
of left translation with the polar map from $K \times \mathfrak k$
to $K^{\mathbb C}$ is a diffeomorphism, and the resulting complex
structure on $\mathrm T^*K\cong \mathrm TK$ combines with the
cotangent bundle symplectic structure to a K\"ahler structure. The
action of $K^{\mathbb C}$ on itself
 by conjugation is holomorphic, and the restriction
of the action to $K$ is Hamiltonian, with momentum mapping from
$K^{\mathbb C}$ to $\mathfrak k^*$ which, viewed as a map on
$K\times \mathfrak k\cong\mathrm T^*K $ and with values in $\mathfrak k$,
amounts to the map
\[
\mu \colon K\times \mathfrak k \longrightarrow \mathfrak k, \quad
\mu(x,Y) = \mathrm{Ad}_xY - Y.
\]
By Proposition 4.2 of \cite{kaehler}, the zero momentum reduced
space
\[
(\mathrm T^*K)_0\cong \mu^{-1}(0)\big/ K
\]
inherits a \emph{stratified K\"ahler structure\/} which is
actually complex algebraic in a sense explained below. Here the
complex algebraic structure is that of the complex algebraic
categorical quotient $K^{\mathbb C}\big/\big/ K^{\mathbb C}$
relative to conjugation, and this quotient, in turn, may be
described as the ordinary orbit space $T^{\mathbb C}\big/W$ of the
complexification $T^{\mathbb C}$ of a maximal torus $T$ of $K$,
relative to the action of the Weyl group $W$; in the literature,
such an orbit space is referred to as an \emph{adjoint
quotient\/}. Certain aspects of the stratified K\"ahler structure
on such an adjoint quotient have been explored in \cite{adjoint}.
The purpose of the present paper is to complement the results in
\cite{adjoint}: For the special case where $K= \mathrm{SU}(n)$ ($n
\geq 2$) we shall elucidate the complex algebraic stratified
K\"ahler structure in terms of complex discriminant varieties and
variants thereof. We will also interpret the resulting stratified
K\"ahler geometry on the adjoint quotient for the special case
where $K=\mathrm{SU}(2)$ in terms of the reduced phase space of a
spherical pendulum constrained to move with angular momentum zero.
In physics, a space of the kind $(\mathrm T^*K)_0$ is the
\emph{building block\/} for certain \emph{lattice gauge\/}
theories. For intelligibility we will include an exposition of the
notion of stratified K\"ahler space.

\section{Stratified K\"ahler spaces}

To develop K\"ahler quantization in the presence of singularities,
we introduced certain \lq\lq K\"ahler spaces with
singularities\rq\rq\ \cite{kaehler} which we refer to as
\emph{stratified K\"ahler spaces\/}. In \cite{qr}, we have shown
that ordinary \emph{K\"ahler quantization\/} can indeed be
extended to a \emph{quantization scheme over stratified K\"ahler
spaces\/}. A special case of a stratified K\"ahler space is a
complex analytic stratified K\"ahler space and, for the present
paper, this notion suffices; we will now describe it.

Let $N$ be a stratified space, the strata being ordinary smooth
manifolds. A \emph{stratified symplectic structure\/} on $N$
consists of a family of symplectic structures, one on each
stratum, together with a Poisson algebra
$(C^{\infty}N,\{\,\cdot\,,\,\cdot\,\})$ of continuous functions on
$N$, and these are required to satisfy the following compatibility
requirements:
\\
{\rm (1)} For each stratum, the restriction map from $C^{\infty}N$
to the algebra of continuous functions on that stratum goes into
the algebra of ordinary smooth functions on the stratum.
\\
{\rm (2)} For each stratum, the restriction map from $C^{\infty}N$
to the algebra of smooth functions on that stratum is a morphism
of Poisson algebras, where the stratum is endowed with the
ordinary smooth symplectic Poisson structure.

A \emph{stratified symplectic space\/} is defined to be a
stratified space together with a stratified symplectic structure.
Given a stratified symplectic space $(N,C^{\infty}N,\{\,\cdot\, ,
\, \cdot\,\})$, the functions in the structure algebra
$C^{\infty}N$ are not necessarily ordinary smooth functions. A
stratified symplectic structure on a space is much more than just
a space stratified into symplectic manifolds: The globally defined
Poisson algebra encapsulates the {\it mutual positions\/} of the
symplectic structures on the strata; in other words, it encodes
the behaviour of the symplectic structures {\it across\/} the
strata.

Recall that a \emph{complex analytic space\/}
(in the sense of {\scsc Grauert\/}) is a topological space $X$, together
with a sheaf of rings $\mathcal O_X$, having the following property:
The space $X$ can be covered by open sets $Y$, each of which
embeds into the open polydisc
$U=\{\mathbf z=(z_1,\dots,z_n);|\mathbf z|<1\}$ in some $\mathbb C^n$
(the dimension $n$ may vary as $U$ varies)
as the zero set of a finite system of holomorphic functions
$f_1,\dots,f_q$ defined on $U$,
such that the
restriction $\mathcal O_Y$ of the sheaf $\mathcal O_X$ to $Y$
is isomorphic as a sheaf to
the quotient sheaf $\mathcal O_U\big/(f_1,\dots,f_q)$;
here $\mathcal O_U$ is the sheaf of germs of holomorphic functions
on $U$. The sheaf $\mathcal O_X$ is then referred to as the
\emph{sheaf of holomorphic functions on \/} $X$.
See \cite{gunnross} for a development of the general theory of complex
analytic spaces.

\begin{defi}
\label{def1}
 A \emph{complex analytic stratified K\"ahler
structure\/} on the stratified space $N$ consists of
\\
{\rm (i)} a stratified symplectic structure
$(C^{\infty}N,\{\,\cdot\,,\,\cdot\,\})$ having the given
stratification of $N$ as its underlying stratification, together
with
\\
{\rm (ii)} a complex analytic structure on $N$ which is compatible
with the stratified symplectic structure.
\end{defi}

Here, the complex analytic structure on $N$ being
\emph{compatible\/} with the stratified symplectic structure means
that the following requirements are met:
\\
(iii) Each stratum is a complex analytic subspace, and the complex
analytic structure, restricted to the stratum, turns that stratum
into an ordinary complex manifold; in particular, the
stratification of $N$ is a refinement of the complex analytic
stratification.
\\
(iv) For each point $q$ of $N$ and each holomorphic function $f$
defined on an open neighborhood $U$ of $q$, there is an open
neighborhood $V$ of $q$ with $V \subset U$ such that, on $V$, $f$
is the restriction of a function in $C^{\infty}(N,\mathbb C) =
C^{\infty}(N) \otimes \mathbb C$.
\\
(v) On each stratum, the symplectic structure  combines with the
complex analytic structure to a K\"ahler structure.

A stratified space $N$, together with a complex analytic
stratified K\"ahler structure, will be said to be a \emph{complex
analytic stratified K\"ahler space\/}.

A simple example of a complex analytic stratified K\"ahler space
arises as follows: In $\mathbb R^3$ with coordinates $x,y,r$,
consider the semicone $N$ given by the equation $x^2 + y^2 = r^2$
and the inequality $r \geq 0$. We refer to this semicone as the
\emph{exotic\/} plane with a single vertex. Consider the algebra
$C^{\infty}N$ of smooth functions in the variables $x,y,r$ subject
to the relation $x^2 + y^2 = r^2$. Thus $C^{\infty}N$ is the
algebra of continuous functions on $N$ which are restrictions of
ordinary smooth functions on the ambient copy of $\mathbb R^3$.
Define the Poisson bracket $\{\,\cdot\,,\,\cdot\,\}$ on this
algebra by
\[
\{x,y\} = 2r,\  \{x,r\} = 2y,\ \{y,r\} = -2x,
\]
and endow $N$ with the complex structure having $z=x+iy$ as
holomorphic coordinate.  The Poisson
and complex analytic structures
combine to a complex analytic stratified K\"ahler structure.
Here the radius function $r$ is \emph{not\/} an ordinary smooth function
of the variables $x$ and $y$. Thus the stratified symplectic
structure cannot be given in terms of ordinary smooth functions of
the variables $x$ and $y$.
The Poisson bracket is \emph{defined
at the vertex\/} as well, away from the vertex the Poisson
structure is an ordinary smooth \emph{symplectic\/} Poisson structure, and
the complex structure does \emph{not\/} \lq\lq see\rq\rq\ the
vertex. Thus the vertex is a
singular point for the Poisson structure whereas it is \emph{not\/}
a singular point for the complex analytic structure.
This semicone $N$ is the classical reduced \emph{phase
space\/} of a single particle moving in ordinary affine space of
dimension $\geq 2$ with angular momentum zero \cite{kaehler},
\cite{varna}.

This example generalizes to an entire class of examples: The
\emph{closure of a holomorphic nilpotent orbit\/} (in a hermitian
Lie algebra) inherits a complex analytic stratified K\"ahler
structure \cite{kaehler}. \emph{Projectivization\/} of  the
closure of a holomorphic nilpotent orbit yields what we refer to
as an \emph{exotic projective variety\/}. In physics, spaces of
this kind arise as reduced classical phase spaces for systems of
harmonic oscillators with zero angular momentum and constant
energy. More details may be found in
\cite{kaehler}--\cite{scorza}, \cite{varna}.

Another class of examples arises from moduli spaces of semistable
holomorphic vector bundles or, more generally, from moduli spaces
of semistable principal bundles on a non-singular complex
projective curve.
See \cite{locpois}--\cite{kaehler} and the literature
there.

Any ordinary K\"ahler manifold is plainly a complex analytic
stratified K\"ahler space. More generally, K\"ahler reduction,
applied to an ordinary K\"ahler manifold, yields  a complex
analytic stratified K\"ahler structure on the reduced space. See
\cite{kaehler} for details. Thus examples of stratified K\"ahler
spaces abound. In the rest of the paper we will explore a
particular class of examples which are actually \emph{algebraic\/}
in a sense which we now explain:

\begin{defi}
\label{def2} Given a stratified space $N$, a \emph{complex
algebraic stratified K\"ahler structure\/} on  $N$ consists of
\\
{\rm (i)} a real semialgebraic structure on $N$ such  that each
stratum inherits the  structure of a real algebraic manifold;
\\
{\rm (ii)} a real algebraic Poisson structure
$\{\,\cdot\,,\,\cdot\,\}$ on (the real structure sheaf of) $N$,
together with a real algebraic symplectic structure on each
stratum, such that the restriction map from (the sheaf of germs of
real algebraic functions on) $N$ to (the sheaf of germs of real
algebraic functions on) each stratum is a Poisson map;
\\
{\rm (iii)} a complex algebraic structure on $N$ which is
compatible with the other structure.
\end{defi}

Given an affine real semialgebraic space $N$, we write its real coordinate
ring as $\mathbb R[N]$. Likewise we write the complex coordinate
ring of an affine complex algebraic variety $N$ as $\mathbb C[N]$.
In this paper, all real semialgebraic structures and all complex algebraic
structures will be affine, and we shall not mention sheaves any more.

At the risk of making a mountain out of a molehill we note that,
in the definition, when the real structure and  the complex
structure on $N$ are both affine (beware: this does not mean that
$\mathbb C[N]$ is the complexification of $\mathbb R[N]$),
the complex algebraic structure on $N$ being
\emph{compatible\/} with the other structure amounts to the
following requirements:
\\
(iv) Each stratum of $N$ is a complex algebraic subspace of $N$,
and the complex algebraic structure, restricted to the stratum,
turns that stratum into an ordinary complex algebraic manifold; in
particular, the stratification of $N$ is a refinement of the
complex algebraic stratification.
\\
(v) For each point $q$ of $N$ and each complex algebraic function
$f$ defined on an open neighborhood $U$ of $q$, there is an open
neighborhood $V$ of $q$ with $V \subset U$ such that, on $V$, $f$
is the restriction of a function in $\mathbb R[N]_{\mathbb C} =
\mathbb R[N] \otimes \mathbb C$.
\\
(vi) On each stratum, the symplectic structure  combines with the
complex algebraic structure to a K\"ahler structure.

A \emph{complex algebraic stratified K\"ahler space\/} is, then, a
stratified space $N$ together with a complex algebraic stratified
K\"ahler structure; the Poisson algebra will then be referred to
as an \emph{algebraic stratified symplectic Poisson algebra\/} (on
$N$).

\section{Complex algebraic stratified K\"ahler structures on adjoint quotients}
\label{castrat}

Return to the situation at the beginning of the paper. We begin
with explaining the K\"ahler structure on $\mathrm T^*K$. Endow
$K$ with the bi-invariant Riemannian metric induced by the
invariant inner product on the Lie algebra $\mathfrak k$. Using
this metric, we identify $\mathfrak k$ with its dual $\mathfrak
k^*$ and the total space of the tangent bundle $\mathrm T K$ with
the total space of the cotangent bundle $\mathrm T^* K$. Thus the
composite
\[
\mathrm T^* K \longrightarrow K \times \mathfrak k \longrightarrow
K^{\mathbb C}
\]
of the inverse of left trivialization with the polar decomposition
map, which assigns $x\cdot\mathrm{exp}(iY) \in K^{\mathbb C}$ to
$(x,Y) \in K \times \mathfrak k$, identifies $\mathrm T^* K$ with
$K^{\mathbb C}$ in a $(K\times K)$-equivariant fashion. Then the
induced complex structure on  $\mathrm T^* K$ combines with the
symplectic structure to a (positive) K\"ahler structure. Indeed,
the {\it real analytic\/} function
\[
\kappa \colon K^{\mathbb C} \longrightarrow \mathbb R, \quad
\kappa(x\cdot\mathrm{exp}(iY)) =  |Y|^2,\quad (x,Y) \in K \times
\mathfrak k,
\]
on $K^{\mathbb C}$ which is twice the \emph{kinetic energy\/}
associated with the Riemannian metric, is a (globally defined)
\emph{K\"ahler potential\/}; in other words, the function $\kappa$
is strictly plurisubharmonic and (the negative of the imaginary
part of) its Levi form yields (what corresponds to) the cotangent
bundle symplectic structure, that is, the cotangent bundle
symplectic structure on $\mathrm T^*K$ is given by
\[
i \partial \overline\partial \kappa=-d \vartheta
\]
where $\vartheta$ is the tautological 1-form on $\mathrm T^*K$. An
explicit calculation which establishes this fact may be found in
\cite{bhallone} (but presumably it is a folk-lore observation). We
note that, given $Y_0 \in \mathfrak k$, the assignment to $(x,Y)
\in K \times \mathfrak k$ of $|Y-Y_0|^2$ yields a K\"ahler
potential as well which, in turn, determines the same K\"ahler
structure as $\kappa$. There is now a literature on related
questions
 \cite{guisteno}, \cite{lempszoe},
\cite{szoekone}.

By Proposition 4.2 in \cite{kaehler}, the zero momentum reduced
space $(\mathrm T^*K)_0$ inherits a complex analytic stratified
K\"ahler structure which is actually a complex algebraic
stratified K\"ahler structure. It is straightforward to describe
this structure directly, and we will now do so: The underlying
complex algebraic structure is that of the categorical quotient $
K^{\mathbb C}\big/\big/ K^{\mathbb C}$ in the category of complex
algebraic varieties, cf. e.~g. \cite{gwschwat}\ (\S 3) for details
on the construction of a categorical quotient. Suffice it to
mention at this stage that $ K^{\mathbb C}\big/\big/ K^{\mathbb
C}$ is the complex affine algebraic variety whose coordinate ring
is the algebra $\mathbb C[K^{\mathbb C}]^{K^{\mathbb C}}$ of
$K^{\mathbb C}$-invariants (relative to the conjugation action) of
the complex affine coordinate ring $\mathbb C[K^{\mathbb C}]$ of
$K^{\mathbb C}$, where  $K^{\mathbb C}$ is viewed as a
(non-singular) complex affine variety. In view of results of
{\scsc Luna\/} \cite{lunatwo}, \cite{lunathr}, the quotient in the
category of algebraic varieties is  the categorical quotient in
the category of analytic varieties as well; for these matters see
also \cite{gwschwat}\ (Theorem 3.6).

Under the present circumstances, the categorical quotient $
K^{\mathbb C}\big/\big/ K^{\mathbb C}$ has a very simple
description: Choose a maximal torus $T$ in $K$ and let $W$ be the
corresponding Weyl group; then $T^{\mathbb C}$ is a maximal torus
in $K^{\mathbb C}$, and the (algebraic) {\it adjoint quotient\/}
$\chi \colon K^{\mathbb C} \to T^{\mathbb C}\big / W$, cf.
\cite{humphone}\ (3.4) and \cite{slodoboo}\ (3.2) for this
terminology, realizes the categorical quotient. By an abuse of
language, we will also refer to the target $T^{\mathbb C}\big / W$
of $\chi$ as \emph{adjoint quotient\/} or as the \emph{adjoint
quotient of\/} $K^{\mathbb C}$. In concrete terms, the map $\chi$
admits the following description: The closure of the conjugacy
class of $x \in K^{\mathbb C}$ contains a unique semisimple
(equivalently: closed) conjugacy class $C_x$ (say), and semisimple
conjugacy classes are parametrized by $T^{\mathbb C}\big / W$; the
image of $x \in K^{\mathbb C}$ under $\chi$ is simply the
parameter value in $T^{\mathbb C}\big / W$ of the semisimple
conjugacy class $C_x$. Since $W$ is a finite group, as a complex
algebraic  space, the quotient $T^{\mathbb C}\big / W$ is simply
the space of $W$-orbits in $T^{\mathbb C}$.

The choice of maximal torus $T$ in $K$ also provides a direct
description of the algebraic stratified symplectic Poisson algebra
on the symplectic quotient $(\mathrm T^* K)_0=\mu^{-1}(0)\big /K$.
Indeed, via the identification (1.1) for $K=T$, the real space
which underlies the (complex algebraic) orbit space $T^{\mathbb
C}\big/W$ for the action of the Weyl group $W$ on $T^{\mathbb C}$
amounts simply to the orbit space $\mathrm T^* T\big/ W$, relative
to the induced action of the Weyl group $W$ on $\mathrm T^* T$,
and the orbit space $\mathrm T^* T\big/ W$ inherits a stratified
symplectic structure in an obvious fashion: Strata are the
$W$-orbits, the closures of the strata are affine varieties, in
the real category as well as in the complex category, indeed,
these closures inherit complex algebraic stratified K\"ahler
structures, and the requisite algebraic stratified symplectic
Poisson algebra $(\mathbb R[\mathrm T^* T\big/W], \{\,\cdot \, ,\,
\cdot\, \})$ is simply the algebra $\mathbb R[\mathrm T^* T]^W$ of
 $W$-invariant functions in $\mathbb R[\mathrm T^* T]$, endowed with the
Poisson bracket $\{\,\cdot \, ,\, \cdot\, \}$ coming from the
ordinary (algebraic) symplectic Poisson bracket on $\mathrm T^*
T$. The choice of invariant inner product on $\mathfrak k$
determines an injection $\mathrm T^* T \to \mathrm T^*K $ which
induces a homeomorphism from $\mathrm T^* T\big/W$ onto $(\mathrm
T^*K)_0$ compatible with all the structure.

\section{The space of normalized degree $n$ polynomials}

In the next section we will work  out defining equations for the
closures of the strata, viewed as complex algebraic varieties, of
the adjoint quotient of $\mathrm{SL}(n,\mathbb C)$. To this end,
we need some preparations, which we now give.

The space $\mathbb A^n_{\mathrm{coef}}$ of normalized degree $n$
polynomials
\begin{equation}
\label{poly} P(z) = z^n + a_1 z^{n-1} + \dots + a_{n-1} z + a_n
\end{equation}
with complex coefficients $a_j$ ($1 \leq j \leq n$) is an
$n$-dimensional complex affine space in an obvious way; it may be
viewed as the $n$'th symmetric power $\mathrm S^n[\mathbb C]$ of a
copy $\mathbb C$ of the complex numbers. The space $\mathbb
A^n_{\mathrm{coef}}$ is stratified according to the multiplicities
of the roots of the polynomials where strata correspond to
partitions of $n$. We define a \emph{partition\/} of $n$ of
\emph{length\/} $r$ to be an $r$-tuple
\begin{equation}
\nu=(n_1, n_2, \ldots ,n_r)
\label{partition}
\end{equation}
of positive natural numbers $n_j$ ($1 \leq j \leq r$) such that
$\sum n_j = n$, normalized so that
\[
n_1 \geq n_2 \geq \ldots\geq n_r.
\]
Given the partition $\nu$ of $n$ of length $r$,
of the kind \eqref{partition},
as a complex manifold,
the stratum corresponding to the partition $\nu$
is
the complex $r$-dimensional manifold $\mathcal D^o_{\nu}$ of polynomials
\[
P(z) = \Pi_{j=1}^r (z-u_j)^{n_j}
\]
with all the roots $u_j$ being pairwise distinct. In particular,
the \emph{top\/} stratum of the space $\mathbb
A^n_{\mathrm{coef}}\cong \mathrm S^n[\mathbb C]$ is the space of
all polynomials having only single roots and, given a length $r$
partition $\nu$ of $n$ of the kind \eqref{partition}, as a complex
manifold, $\mathcal D^o_{\nu}$ comes down to the top stratum of
$\mathrm S^r[\mathbb C]$, that is, to the space of normalized
degree $r$ polynomials having only single roots. For each
partition $\nu$ of $n$ of the kind \eqref{partition}, let
$\mathcal D_{\nu}$ be the closure of $\mathcal D^o_{\nu}$ in
$\mathbb A^n_{\mathrm{coef}}$; this closure is an affine variety.
In particular, at the bottom, we have $\mathcal D_{(n)}=\mathcal
D^o_{(n)}$. Given the partition $\nu=(n_1, n_2, \ldots ,n_r)$ of
$n$, consider the partition of $n$ which arises from $\nu$ by the
operation of picking a pair $(n_j,n_k)$, taking the sum of $n_j$
and $n_k$, and reordering the terms (if need be); iterating this
operation we obtain a partial ordering among the partitions of
$n$, which we write as $\nu' \preceq \nu$, where $\nu'$ arises
from $\nu$ by a finite sequence of operations of the kind just
described, the empty sequence being admitted, so that $\nu \preceq
\nu$. When $\nu' \preceq \nu$ and $\nu' \neq \nu$, we write $\nu'
\prec \nu$. The stratum $\mathcal D^o_{\nu'}$ lies in the closure
$\mathcal D_{\nu}$ of $\mathcal D^o_{\nu}$ if and only if $\nu'
\prec \nu$; thus, for any partition $\nu$ of $n$,
\[
\mathcal D_{\nu} = \cup_{\nu' \preceq \nu}\mathcal D^o_{\nu'}.
\]
It is manifest that $\mathbb A^n_{\mathrm{coef}}$
is the disjoint union of the $\mathcal D^o_{\nu}$'s
as $\nu$ ranges over partitions of $n$, and this
decomposition is a stratification.

The stratification of $\mathbb A^n_{\mathrm{coef}}$
can be understood in terms of
discriminant varieties.
We will now explain this briefly;
 cf. \cite{gkatzone} and the literature there for more details.

Let $\Sigma_1$ be the non-singular hypersurface in the
$(n+1)$-dimensional complex affine space $\mathbb C \times \mathbb
A^n_{\mathrm{coef}}$ with coordinates $(z,a_1,\ldots,a_{n-1},a_n)$
given by the equation
\begin{equation}
 z^n + a_1 z^{n-1} + \dots + a_{n-1} z + a_n = 0.
\end{equation}
Then $a_n$ is a polynomial function of the other coordinates
$z,a_1,\ldots,a_{n-1}$ whence the hypersurface $\Sigma_1$ admits
the parametrization
\begin{equation}
\label{para}
h \colon \mathbb C^n \to
\mathbb C \times \mathbb A^n_{\mathrm{coef}},
\quad
(z,a_1,\ldots,a_{n-1}) \mapsto (z,a_1,\ldots,a_{n-1},a_n),
\end{equation}
and this parametrization is non-singular (smooth).

Given the normalized complex polynomial $P(z)$ of the kind
\eqref{poly}, for $1 \leq k < n$, consider the successive
derivatives
\[
P^{(k)}(z) = k! z^{n-k} + (k-1)!a_1 z^{n-k-1} + \dots + a_{n-k};
\]
accordingly, let $\Sigma_{k+1} \subseteq \Sigma_1$ be the affine
complex variety given by the equations
\begin{equation}
P(z) =0,\ P'(z) =0,\ \ldots,\  P^{(k)}(z) =0;
\label{equations}
\end{equation}
the variety $\Sigma_{k+1}$ is actually non-singular. Let
\[
p \colon \mathbb C \times \mathbb A^n_{\mathrm{coef}}
\longrightarrow
\mathbb A^n_{\mathrm{coef}}
\]
be the projection to the second factor and, for $1 \leq k \leq n$,
let $\mathcal D_k = p(\Sigma_k)\subseteq \mathbb
A^n_{\mathrm{coef}}$, so that, in particular, $\mathcal
D_1=\mathbb A^n_{\mathrm{coef}}$; we refer to $\mathcal D_k$ as
the $k$'th \emph{discriminant\/} variety. We shall see shortly
that each $\mathcal D_k$ is indeed a complex algebraic variety.
For $1 \leq k \leq n-1$, the restriction of the projection $p$ to
$\Sigma_{k} \setminus \Sigma_{k+1}$ is $(n-k+1)$ to $1$, the
restriction of the projection $p$ from $\Sigma_{k}$ to $\mathcal
D_k$ branches over $\mathcal D_{k+1} \subseteq \mathcal D_k$, and
the restriction of the projection $p$ to $\Sigma_{n}$ is injective
and identifies $\Sigma_{n}$ with $\mathcal D_n$. For $1 \leq k
\leq n$, $\mathcal D_k$ is the space of polynomials having at
least one root with multiplicity at least equal to $k$. In
particular, $\mathcal D_n$ is a smooth curve, parametrized by
\[
w \mapsto\left(-nw, \binom n 2 w^2, -\binom n3 w^3, \ldots, (-1)^n w^n
\right) ,\quad w \in \mathbb C.
\]
By construction, the $\mathcal D_k$'s form an ascending sequence
\[
\mathcal D_n \subseteq
\mathcal D_{n-1} \subseteq \ldots \subseteq
\mathcal D_1 =\mathbb A^n_{\mathrm{coef}}.
\]

Equations for the $\mathcal D_k$'s can be obtained by the
following classical procedure: Given the polynomial
\begin{equation}
\label{poly2}
P(z) = a_0z^n + a_1 z^{n-1} + \dots + a_{n-1} z + a_n,
\quad a_j \in \mathbb C \quad (1 \leq j \leq n),
\end{equation}
with roots $z_1,\dots,z_n$,
the \emph{discriminant\/}
\[
D_n(P)= D_n(a_0,a_1,\dots,a_n)
\]
of this polynomial is defined to be the expression
\begin{equation}
\label{discri}
D_n(a_0,a_1,\dots,a_n)
=
a_0^{2n-2}
\Pi_{j<k}(z_j -z_k)^2,
\end{equation}
written as a polynomial in
the coefficients
$a_0,a_1,\dots,a_n$.

To reproduce an expression
for the discriminant $D_n(a_0,a_1,\dots,a_n)$
of a polynomial of the kind \eqref{poly2},
recall from e.~g. \cite{vandewae} that, given two polynomials
\begin{align}
f(x) &= a_0 x^n + a_1 x^{n-1} + \ldots + a_n
\\
h(x) &= b_0 x^m + b_1 x^{m-1} + \ldots + b_m,
\end{align}
the \emph{resultant\/} $R(f,h)$ of $f$ and $h$, sometimes referred
to as the {\scsc Sylvester\/} \emph{resultant\/},
 may be computed as
the $((m+n)\times(m+n))$-determinant
\begin{equation}
R(f,h)=\left|\begin{array}{cccccccccc}
a_0  & a_1 & a_2 &\cdots&  a_n  &  0  & 0     &\cdots& 0    &0
\\
0    & a_0 & a_1 &\cdots&a_{n-1}& a_n & 0     &\cdots& 0    &0
\\
\cdot&\cdot&\cdot&\cdots& \cdot &\cdot& \cdot &\cdots&\cdot &\cdot
\\
\cdot&\cdot&\cdot&\cdots& \cdot &\cdot& \cdot &\cdots&\cdot &\cdot
\\
\cdot&\cdot&\cdot&\cdots& \cdot &\cdot& \cdot &\cdots&\cdot &\cdot
\\
0    &  0  &  0  &  0   &  a_0  & a_1 & a_2   &\cdots&a_{n-1}& a_n
\\
b_0  & b_1 & b_2 &\cdots&  \cdot& b_m  &  0   & 0    &\cdots & 0
\\
0    & b_0 & b_1 &\cdots&\cdot  &b_{m-1}& b_m & 0    &\cdots & 0
\\
\cdot&\cdot&\cdot&\cdots& \cdot &\cdot  &\cdot&\cdots&\cdot&\cdot
\\
\cdot&\cdot&\cdot&\cdots& \cdot &\cdot& \cdot &\cdots&\cdot &\cdot
\\
\cdot&\cdot&\cdot&\cdots& \cdot &\cdot& \cdot &\cdots&\cdot &\cdot
\\
0    &  0  &  0  &  0   &  b_0  & b_1   & b_2 &\cdots&b_{m-1}& b_m
\end{array}
\right|.
\end{equation}
This determinant has $m$ rows of the kind
\[
[0,\ldots, 0,a_0, a_1, a_2,\cdots, a_n,0,\ldots, 0]
\]
and $n$ rows of the kind
\[
[0,\ldots, 0,b_0, b_1, b_2,\cdots, b_m,0,\ldots, 0]
\]
We quote the following classical facts \cite{vandewae}.

\begin{prop}
When $a_0 \ne 0 \ne b_0$, the polynomials $f$ and $h$ have a
common linear factor if and only if $R(f,h) = 0$.
\end{prop}

\begin{prop} The discriminant $D_n(a_0,\ldots,a_n)$ of the polynomial
\[
f(x) = a_0 x^n + a_1 x^{n-1} + \ldots + a_n
\]
 satisfies the identity
\begin{equation}
a_0D = \pm R(f,f').
\end{equation}
\end{prop}

Now, by construction, for $1 \leq k < n$,
the discriminant variety $\mathcal D_{k+1} \subseteq \mathcal D_1$
is given by the equations
\begin{equation}
\label{discri2}
\begin{aligned}
D_n(1,a_1,\ldots,a_n) &=0
\\
D_{n-1}(n ,(n-1)a_1,\ldots,a_{n-1}) &=0
\\
{}&\ldots
\\
D_{n-{(k-1)}}\left(\frac {n!}{(n-k)!},
\frac{(n-1)!}{(n-k)!}a_1,\ldots,a_{n-{(k-1)}}\right) &=0 ;
\end{aligned}
\end{equation}
the dimension of
$\mathcal D_{k+1}$ equals $n-k$.

For illustration, we give explicit expressions for
the descriminants for $2 \leq n \leq 4$:
When $n=2$, the discriminant $D_2(a_0,a_1,a_2)$ of
$a_0x^2 + a_1x + a_2$ comes
down to the familiar expression
\begin{equation}
D_2(a_0,a_1,a_2)=a_1^2 -4 a_0a_2;
\end{equation}
and when $n=3$, the discriminant
$D_3(a_0,a_1,a_2,a_3)$ of the cubic polynomial
$a_0x^3 + a_1x^2 + a_2x + a_3$
reads
\begin{equation}
D_3(a_0,a_1,a_2,a_3) =a_1^2a_2^2 -4 a_0a_2^3 -4 a_1^3 a_3 -27
a_0^2a_3^2 + 18 a_0 a_1 a_2 a_3 . \label{cubic}
\end{equation}
Likewise, the discriminant $D_4=D_4(a_0,a_1,a_2,a_3,a_4)$ of
the quartic polynomial
$a_0x^4 + a_1x^3 + a_2x^2 + a_3x + a_4$
has the form
\[
\begin{aligned}
D_4&= (a_1^2 a_2^2 a_3^2 - 4 a_1^3 a_3^3 -4 a_0 a_2^3 a_3^2
+18a_0 a_1 a_2 a_3^3
-27 a_0^2 a_3^4 + 256 a_0^3 a_4^3)
\\
& \quad + (-4a_1^2 a_2^3 + 18 a_1^3 a_2 a_3 + 16 a_0 a_2^4
\\
& \quad \quad  - 80 a_0 a_1 a_2^2 a_3 -6 a_0 a_1^2 a_3^2 + 144
a_0^2 a_2 a_3^2)a_4
\\
& \quad
+(-27 a_1^4 + 144 a_0 a_1^2 a_2 - 128 a_0^2 a_2^2 - 192 a_0^2 a_1 a_3)a_4^2 .
\end{aligned}
\]

For $1 \leq k \leq n$,
when $\nu$ is the partition
$(k,1,\ldots,1)$ of $n$ of length $k-r+1$,
$\mathcal D_{\nu}$ equals $\mathcal D_{k}$;
moreover, the closure $\mathcal D_{k}$ of $\mathcal D^o_{\nu}$
then contains the $\mathcal D^o_{\nu}$'s for
$\nu=(n_1, n_2,\ldots, n_s\}$ with $n_1 \geq k$.
Given $\mathcal D_{\nu}$ with
$\nu=\{k, n_2, \ldots, n_r\}$,
defining equations for $\mathcal D_{\nu}$
may be obtained by adding suitable equations to those defining
$\mathcal D_{k}$.

In low dimensions, these considerations yields the following picture:

\noindent $n=2$: In this case, the variety
$\mathbb A^2_{\mathrm{coef}}=\mathcal D_1$
is a 2-dimensional
complex affine space,
which decomposes into
two strata according to the two partitions
$(1,1)$ and $(2)$ of the natural number $n=2$.
In the coordinates $(a_1,a_2)$ on
$\mathbb A^2_{\mathrm{coef}}$,
the variety $\mathcal D_2$ is the non-singular curve of degree 2
given by the equation
\begin{equation}
D_2(1,a_1,a_2) = a_1^2 - 4 a_2 = 0 .
\label{equation3}
\end{equation}

\noindent $n=3$: Now the variety $\mathbb
A^3_{\mathrm{coef}}=\mathcal D_1$ is a 3-dimensional complex
affine space, which decomposes into three strata according to the
three partitions $(1,1,1)$, $(2,1)$, $(3)$ of $n=3$. Furthermore,
in the coordinates $(a_1,a_2,a_3)$ on $\mathbb
A^3_{\mathrm{coef}}$, the surface $\mathcal D_2= \mathcal
D_{(2,1)}$ in $\mathbb A^3_{\mathrm{coef}}$ is given by the
equation
\begin{equation}
D_3(1,a_1,a_2,a_3) =  0 ,
\label{equation4}
\end{equation}
and this surface
decomposes as
\[
\mathcal D_2 = \mathcal D^o_{(2,1)} \cup
\mathcal D_{3}.
\]
Moreover, in the chosen coordinates,
$\mathcal D_3$ is the curve in $\mathbb A^3_{\mathrm{coef}}$
given by the equations
\eqref{equation4} and
\begin{equation}
D_2(3,2a_1,a_2) =  0 ,
\label{equation5}
\end{equation}
and this curve is non-singular. The surface $\mathcal D_2 =
\mathcal D_{(2,1)}$ in $\mathbb A^3_{\mathrm{coef}}$ is singular
along the curve $\mathcal D_3$. Indeed, a calculation shows that
the three partial derivatives of the defining equation
\eqref{equation4} of $\mathcal D_2$ vanish identically along
$\mathcal D_3$. A closer look reveals that $\mathcal D_2$ has a
\lq\lq fold\rq\rq\  along $\mathcal D_3$ but is topologically flat
in $\mathbb A^3_{\mathrm{coef}}$; cf. Corollary 5.1 of
\cite{gkatzone} for details.

\noindent $n=4$: In this case, by construction, $\mathbb
A^4_{\mathrm{coef}}=\mathcal D_1$ is a 4-dimensional complex
affine space; this space decomposes into five strata according to
the five partitions $(1,1,1,1)$, $(2,1,1)$, $(2,2)$, $(3,1)$,
$(4)$ of $n=4$. Furthermore, in the coordinates
$(a_1,a_2,a_3,a_4)$ on $\mathbb A^4_{\mathrm{coef}}$, the
hypersurface $\mathcal D_2= \mathcal D_{(2,1,1)}$ in $\mathcal
D_1$ is given by the equation
\begin{equation}
D_4(1,a_1,a_2,a_3,a_4) =  0 ,
\label{equation6}
\end{equation}
and this hypersurface
decomposes as
\[
\mathcal D_{(2,1,1)} =\mathcal D^o_{(2,1,1)}
\cup \mathcal D^o_{(2,2)} \cup \mathcal D^o_{(3,1)}
\cup  \mathcal D_{4}.
\]
Likewise, in the chosen coordinates,
the surface
$\mathcal D_3= \mathcal D_{(3,1)}$ in
$\mathbb A^4_{\mathrm{coef}}$
is given by the equations \eqref{equation6}
and
\begin{equation}
D_3(4,3a_1,2a_2,a_3) =  0 ,
\label{equation7}
\end{equation}
and this surface decomposes as
\[
\mathcal D_{(3,1)} =\mathcal D^o_{(3,1)} \cup
\mathcal D_{4}.
\]
Moreover,
$\mathcal D_4$ is a non-singular curve of degree 4 in
$\mathbb A^4_{\mathrm{coef}}$
which, in the above coordinates, is
given by the equations
\eqref{equation6}, \eqref{equation7}, and
\begin{equation}
D_2(12,6a_1,2a_2) =  0 .
\label{equation8}
\end{equation}
The hypersurface $\mathcal D_2=\mathcal D_{(2,1,1)}$
in $\mathbb A^4_{\mathrm{coef}}$
is the space of normalized complex degree 4 polynomials
with at least one multiple root, and
the singular locus of this hypersurface $\mathcal D_{(2,1,1)}$
is the union
\[
\mathcal D_3 \cup \mathcal D_{(2,2)}
=\mathcal D^o_{(3,1)} \cup \mathcal D^o_{(2,2)}
\cup \mathcal D_{(4)} .
\]
The surface $\mathcal D_3=\mathcal D_{(3,1)}$
in
$\mathbb A^4_{\mathrm{coef}}$
is the space of normalized complex degree 4 polynomials
with a root of multiplicity at least 3,
and the singular locus of
this surface $\mathcal D_3$
is the curve $\mathcal D_4$.
See Example 6.1 in \cite{gkatzone}
for details.
The surface
\[
\mathcal D_{(2,2)} =\mathcal D^o_{(2,2)} \cup  \mathcal D_{(4)}
\]
in
$\mathbb A^4_{\mathrm{coef}}$
is the space of normalized complex degree 4 polynomials
having two roots with multiplicity 2 or a single root with multiplicity 4;
this surface
is \emph{not\/} a discriminant variety of the kind $\mathcal D_k$.
Inspection
of the coefficients $a_1,a_2,a_3,a_4$
of the polynomial
\[
P(z) = (z-z_1)^2(z-z_2)^2 = z^4 + a_1 z^3 + a_2 z^2 + a_3 z + a_4
\]
shows that,
in the coordinates $(a_1,a_2,a_3,a_4)$,
 the variety
$\mathcal D_{(2,2)}$
in
$\mathbb A^4_{\mathrm{coef}}$
is given by the two equations
\begin{align}
(a_1^2-4a_2) a_1 + 4a_3 &=0
\label{D21}
\\
(a_1^2-4a_2)^2 - 16a_4 &=0 .
\label{D22}
\end{align}
Thus $a_3$ and $a_4$ are polynomial functions of $a_1$ and $a_2$
whence the obvious parametrization
\[
(a_1,a_2) \longmapsto (a_1,a_2,a_3,a_4)
\]
identifies
$\mathcal D_{(2,2)}$ with a copy of 2-dimensional
complex affine space $\mathbb A^2$.
Inspection
of the coefficients $a_1$ and $a_2$
of the polynomial
\begin{equation}
P(z) = (z-z_0)^4 = z^4 + a_1 z^3 + a_2 z^2 + a_3 z + a_4
\label{poly3}
\end{equation}
shows that,
in terms of the coordinates $a_1$ and $a_2$,
the variety $\mathcal D_4\subseteq \mathcal D_{(2,2)}$ is the curve
given by the equation
\begin{equation}
3a_1^2-8a_2 = 0.
\label{equation2}
\end{equation}
We intend to work out elsewhere how equations for varieties of the
kind $\mathcal D_{\nu}$ which are not ordinary discriminant
varieties (i.~e. not of the kind $\mathcal D_{k}$) may be derived.

\section{The adjoint quotient
of $\mathrm{SL}(n,\mathbb C)$} \label{SL}

Let $n \geq 2$. In this section we will explore the stratification
of the adjoint quotient of $\mathrm{SL}(n,\mathbb C)$. For
convience, we will consider first the group $\mathrm{GL}(n,\mathbb
C)$. A maximal complex torus in $\mathrm{GL}(n,\mathbb C)$ is
given by the complex diagonal matrices in $\mathrm{GL}(n,\mathbb
C)$. Thus, as a complex Lie group, this torus  is isomorphic to
the product $(\mathbb C^*)^n$ of $n$ copies of the multiplicative
group $\mathbb C^*$ of non-zero complex numbers. Introduce
standard coordinates $z_1,\dots,z_n$ on $\mathbb C^n$. Then the
fundamental characters $\sigma_1, \dots, \sigma_n$ of
$\mathrm{GL}(n,\mathbb C)$, restricted to the maximal torus, come
down to  the elementary symmetric functions in the variables
$z_1,\dots,z_n$, and the assignment to $\mathbf z
=(z_1,\dots,z_n)$ of
\[
(a_1,\dots,a_n) = (-\sigma_1(\mathbf z),\sigma_2(\mathbf z),\ldots,
(-1)^n\sigma_n(\mathbf z)) \in \mathbb C^n
\]
yields a holomorphic map
\begin{equation}
 (-\sigma_1,\sigma_2\dots,(-1)^n\sigma_n)\colon
(\mathbb C^*)^n
 \longrightarrow
\mathbb C^n
\label{steinberg1}
\end{equation}
which induces a complex algebraic isomorphism from the adjoint
quotient $(\mathbb C^*)^n\big/S_n$ of $\mathrm{GL}(n,\mathbb C)$
onto the subspace of the space $\mathbb A^n_{\mathrm{coef}}$ of
complex normalized degree $n$ polynomials which consists of
normalized polynomials having non-zero constant coefficient.

Let $K=\mathrm{SU}(n)$, so
that $K^{\mathbb C}=\mathrm{SL}(n,\mathbb C)$.
A maximal complex torus $T^{\mathbb C}$
in $\mathrm{SL}(n,\mathbb C)$
is given by the complex diagonal matrices in
$\mathrm{SL}(n,\mathbb C)$, that is,
 by the complex diagonal $(n \times n)$-matrices having determinant 1.
Realize the torus $T^{\mathbb C}$ as the subspace of
the maximal torus of
$\mathrm{GL}(n,\mathbb C)$
in the standard fashion, that is,
as the subspace of $(\mathbb C^*)^n$
which consists of all $(z_1,\dots,z_n)$ in $(\mathbb
C^*)^n$ such that $z_1 \dots z_n = 1$.
The
holomorphic map \eqref{steinberg1}
restricts to the holomorphic
map
\begin{equation}
 (-\sigma_1,\dots,(-1)^{n-1}\sigma_{n-1})\colon
T^{\mathbb C} \cong
(\mathbb C^*)^{n-1}
 \longrightarrow
\mathbb C^{n-1} .
\label{steinberg2}
\end{equation}
This map induces a complex algebraic isomorphism from the adjoint
quotient $T^{\mathbb C}\big/S_n$ of $\mathrm{SL}(n,\mathbb C)$
onto the subspace $\mathbb A^{n,1}_{\mathrm{coef}}$ of the space
$\mathbb A^n_{\mathrm{coef}}$ of complex normalized degree $n$
polynomials which consists of polynomials having constant
coefficient equal to 1. The space $\mathbb
A^{n,1}_{\mathrm{coef}}$ is a complex affine space of dimension
$n-1$ whence the same is true of the adjoint quotient of
$\mathrm{SL}(n,\mathbb C)$. The discussion in the previous section
applies to the adjoint quotient of $\mathrm{SL}(n,\mathbb C)$ in
the following fashion:

For $1 \leq k \leq n$, let $ \mathcal D^{\mathrm {SU}(n)}_k=
\mathcal D_k\cap \mathbb A^{n,1}_{\mathrm{coef}} $ so that, in
particular, $\mathcal D^{\mathrm {SU}(n)}_1$ is the entire space
$\mathbb A^{n,1}_{\mathrm{coef}} \cong \mathrm{SL}(n,\mathbb
C)\big/ \big/ \mathrm{SL}(n,\mathbb C). $ Likewise, for each
partition $\nu$ of $n$,  let
\[
\mathcal D^{\mathrm {SU}(n)}_{\nu}=
\mathcal D_{\nu} \cap
\mathbb A^{n,1}_{\mathrm{coef}}
\]
and
\[
\mathcal D^{(o,\mathrm {SU}(n))}_{\nu}=
\mathcal D^o_{\nu} \cap
\mathbb A^{n,1}_{\mathrm{coef}}.
\]
Then $\mathbb A^{n,1}_{\mathrm{coef}}\cong \mathrm{SL}(n,\mathbb
C)\big/ \big/ \mathrm{SL}(n,\mathbb C)$ is the disjoint union of
the $\mathcal D^{(o,\mathrm {SU}(n))}_{\nu}$'s as $\nu$ ranges
over partitions of $n$, and this decomposition is a
stratification. In low dimensions, we are thus led to the
following picture:

\noindent $n=2$: In this case, the adjoint quotient $\mathbb
A^{2,1}_{\mathrm{coef}}=\mathcal D^{\mathrm{SU}(2)}_1$ is a copy
of the complex line, which decomposes into two strata according to
the two partitions $(1,1)$ and $(2)$ of the natural number $n=2$.
In  the coordinate $a_1$ on $\mathbb A^{2,1}_{\mathrm{coef}}$, the
stratum $\mathcal D^{\mathrm{SU}(2)}_2$ consists of the two points
$\pm 2$; these are solutions of the equation
\begin{equation}
D_2(1,a_1,1) = a_1^2 - 4 = 0 .
\label{equation13}
\end{equation}

\noindent $n=3$: Now the adjoint quotient
$\mathbb A^{3,1}_{\mathrm{coef}}=\mathcal D^{\mathrm{SU}(3)}_1$
is a 2-dimensional
complex affine space, which decomposes into
three strata according to the three partitions
$(1,1,1)$, $(2,1)$, $(3)$ of $n=3$.
Furthermore, in the coordinates $(a_1,a_2)$ on
$\mathbb A^{3,1}_{\mathrm{coef}}$, the variety
$\mathcal D^{\mathrm{SU}(3)}_2= \mathcal D^{\mathrm{SU}(3)}_{(2,1)}$ in
$\mathbb A^{3,1}_{\mathrm{coef}}$
is the curve given by the equation
\begin{equation}
D_3(1,a_1,a_2,1) =  0 ,
\label{equation14}
\end{equation}
that is, in view of \eqref{cubic}, by the equation
\begin{equation}
a_1^2a_2^2 -4 a_2^3 -4 a_1^3 -27 + 18 a_1
a_2 =0,
\end{equation}
and this curve
decomposes as
\[
\mathcal D^{\mathrm{SU}(3)}_{(2,1)}
= \mathcal D^{(\mathrm{SU}(3),o)}_{(2,1)} \cup
\mathcal D^{\mathrm{SU}(3)}_{3}.
\]
Moreover, in the chosen coordinates,
the lowest stratum $\mathcal D^{\mathrm{SU}(3)}_3$
consists of the three points
in $\mathbb A^3_{\mathrm{coef}}$
which solve  the equations
\eqref{equation14} and
\begin{equation}
D_2(3,2a_1,1) =  0 .
\label{equation15}
\end{equation}
At each of these three points, the curve $\mathcal
D^{\mathrm{SU}(3)}_2$ in $\mathbb A^3_{\mathrm{coef}}$ has a cusp,
as a direct examination shows. This fact is also a consequence of
the observation spelled out in the previous section that the
surface $\mathcal D_2$ in $\mathbb A^3_{\mathrm{coef}}$ has a
\lq\lq fold\rq\rq\  along $\mathcal D_3$. Topologically, the curve
$\mathcal D^{\mathrm{SU}(3)}_2$ in $\mathbb A^3_{\mathrm{coef}}$
is flat, i.~e. the embedding of $\mathcal D^{\mathrm{SU}(3)}_2$
into $\mathbb A^3_{\mathrm{coef}}$ is locally homeomorphic to the
standard embedding of $\mathbb R^2$ into $\mathbb R^6$.

\noindent $n=4$: In this case, the adjoint quotient
$\mathbb A^{4,1}_{\mathrm{coef}}=\mathcal D^{\mathrm{SU}(4)}_1$
is a 3-dimensional complex affine space,
which decomposes into five strata according to the five partitions
$(1,1,1,1)$, $(2,1,1)$, $(2,2)$, $(3,1)$, $(4)$ of $n=4$.
Furthermore, in the coordinates $(a_1,a_2,a_3)$ on
$\mathbb A^{4,1}_{\mathrm{coef}}$,
the surface
$\mathcal D^{\mathrm{SU}(4)}_2= \mathcal D^{\mathrm{SU}(4)}_{(2,1,1)}$
in $\mathbb A^{4,1}_{\mathrm{coef}}$ is given by the equation
\begin{equation}
D_4(1,a_1,a_2,a_3,1) =  0 ,
\label{equation16}
\end{equation}
and this surface
decomposes as
\[
\mathcal D^{\mathrm{SU}(4)}_{(2,1,1)}
=\mathcal D^{(\mathrm{SU}(4),o)}_{(2,1,1)}
\cup \mathcal D^{(\mathrm{SU}(4),o)}_{(2,2)}
\cup \mathcal D^{(\mathrm{SU}(4),o)}_{(3,1)}
\cup  \mathcal D^{\mathrm{SU}(4)}_{4}.
\]
Likewise, in the chosen coordinates,
the curve
$\mathcal D^{\mathrm{SU}(4)}_3= \mathcal D^{\mathrm{SU}(4)}_{(3,1)} $ in
$\mathbb A^{4,1}_{\mathrm{coef}}$
is given by the equations \eqref{equation16}
and
\begin{equation}
D_3(4,3a_1,2a_2,a_3) =  0 ,
\label{equation17}
\end{equation}
and this curve decomposes as
\[
\mathcal D^{\mathrm{SU}(4)}_{(3,1)}
=\mathcal D^{(\mathrm{SU}(4),o)}_{(3,1)} \cup
\mathcal D^{\mathrm{SU}(4)}_{4}.
\]
Moreover, the lowest stratum $\mathcal D^{\mathrm{SU}(4)}_4$
consists of the four points in
$\mathbb A^{4,1}_{\mathrm{coef}}$
which, in the coordinates $(a_1,a_2,a_3)$, solve
the equations
\eqref{equation16}, \eqref{equation17}, and
\begin{equation}
D_2(12,6a_1,2a_2) =  0 .
\label{equation18}
\end{equation}
The singular locus of the surface $\mathcal D^{\mathrm{SU}(4)}_2=
 \mathcal D^{\mathrm{SU}(4)}_{(2,1,1)} $
in $\mathbb A^{4,1}_{\mathrm{coef}}$
is the union
\[
\mathcal D^{\mathrm{SU}(4)}_3 \cup \mathcal D^{\mathrm{SU}(4)}_{(2,2)}
=\mathcal D^{(\mathrm{SU}(4),o)}_{(3,1)} \cup
\mathcal D^{(\mathrm{SU}(4),o)}_{(2,2)}
\cup \mathcal D^{\mathrm{SU}(4)}_{4} .
\]
The singular locus of
the curve $\mathcal D^{\mathrm{SU}(4)}_3$ in
$\mathbb A^{4,1}_{\mathrm{coef}}$
consists of the four points in $\mathcal D^{\mathrm{SU}4}_4$.
This is a consequence of the corresponding observation
describing the singular locus
of the surface $\mathcal D_3$ in $\mathbb A^4_{\mathrm{coef}}$,
spelled out in the previous section.
The curve
\[
\mathcal D^{\mathrm{SU}(4)}_{(2,2)}
=\mathcal D^{(\mathrm{SU}(4),o)}_{(2,2)} \cup
\mathcal D^{\mathrm{SU}(4)}_{4}
\]
in
$\mathbb A^{4,1}_{\mathrm{coef}}$
is given by the two equations
\begin{align}
(a_1^2-4a_2) a_1 + 4a_3 &=0
\label{D221}
\\
(a_1^2-4a_2)^2 - 16 &=0 .
\label{D222}
\end{align}
Consequently the embedding
\[
(a_1,a_2) \longmapsto (a_1,a_2,\frac 14 (4a_2- a_1^2) a_1)
\]
of 2-dimensional complex affine space $\mathbb A^2$ into $\mathbb
A^{4,1}_{\mathrm{coef}}$ identifies the non-sin\-gular curve in
$\mathbb A^2$ given by the equation \eqref{D222} with the curve
$\mathcal D^{\mathrm{SU}(4)}_{(2,2)}$ in $\mathbb
A^{4,1}_{\mathrm{coef}}$ which is therefore neccessarily
non-singular. The observation related with the polynomial
\eqref{poly3} in the previous section shows that, in terms of the
coordinates $a_1$ and $a_2$, the bottom stratum $\mathcal
D^{\mathrm{SU}(4)}_4$, viewed as a subset of $\mathcal
D^{\mathrm{SU}(4)}_{(2,2)}$, consists of the four points in
$\mathcal D^{\mathrm{SU}(4)}_{(2,2)}$ solving the equation
\eqref{equation2}, that is, the four points given by $(a_1,a_2) =
(\pm 2 \sqrt 2, \pm 3)$.

It is worth remarking on the position of the real orbit space
$\mathrm {SU}(n)\big/\mathrm {SU}(n)$ relative to conjugation,
realized within the adjoint quotient $T^{\mathbb C}\big/S_n$.
Thus, consider the real affine space $\mathbb R^{n-1}$, embedded
into $\mathbb C^{n-1}$ as the real affine space of real points of
the adjoint quotient $\mathbb C^{n-1}\cong T^{\mathbb C}\big/S_n$
in the obvious fashion. The space
\begin{equation}
\mathcal H_n = \mathcal D^{\mathrm{SU}(n)}_2 \cap \mathbb R^{n-1}
\label{hypersurface}
\end{equation}
of real points of the complex variety $\mathcal
D^{\mathrm{SU}(n)}_2$ is a real compact hypersurface in $\mathbb
R^{n-1}$, and the orbit space $\mathrm {SU}(n)\big/\mathrm
{SU}(n)$, realized within the adjoint quotient $T^{\mathbb
C}\big/S_n\cong C^{n-1}$ of $\mathrm {SL}(n,\mathbb C)$, amounts
to the compact region $R_n$ in $\mathbb R^{n-1}\subseteq \mathbb
C^{n-1}$ bounded by this hypersurface. The exponential mapping
from $\mathfrak t$ to $T$, restricted to a Weyl alcove, induces a
homeomorphism from this alcove onto $R_n$.

\section{The stratified K\"ahler structure on the adjoint quotient
of $\mathrm{SL}(n,\mathbb C)$}

Let $K$ be a general compact Lie group, let $T$ be a maximal torus
in $K$, let $K^{\mathbb C}$ and $T^{\mathbb C}$ be the
complexification of $K$ and $T$, respectively,
 and let $N$ be the adjoint quotient of $K^{\mathbb C}$, that is,
$N$ equals the space $T^{\mathbb C}\big/W$ of $W$-orbits relative
to the action of the Weyl group $W$ of $K$ on $T^{\mathbb C}$. The
algebraic stratified symplectic Poisson algebra on the
\emph{real\/} coordinate ring $\mathbb R[T^{\mathbb C}]$ of the
complex torus $T^{\mathbb C}$ is obtained in the following
fashion: The torus $T^{\mathbb C}$ amounts to a product of
finitely many copies of the complex 1-dimensional torus $\mathbb
C^*$. In terms of the standard coordinate $z= x+ i y$ on $\mathbb
C^*$, identified with the multiplicative group of non-zero complex
numbers, the corresponding real Poisson structure is given by
$\{x,y\} = x^2 + y^2$, and the resulting symplectic structure and
the complex structure on $\mathbb C^*$ combine to an algebraic
K\"ahler structure which, in turn, is just the structure arising
from the standard embedding of $\mathbb C^*$ into
$\mathrm{SU}(2)^{\mathbb C}= \mathrm{SL}(2,\mathbb C)$ as a
maximal complex torus, where $\mathrm{SL}(2,\mathbb C)$ is endowed
with the K\"ahler structure reproduced in Section \ref{castrat}
above. The symplectic Poisson algebra on the complex torus
$T^{\mathbb C}$ in the general group  $K^{\mathbb C}$
is now simply the product structure whence the
induced symplectic structure on  $T^{\mathbb C}$ is plainly real algebraic. The
algebraic stratified symplectic Poisson algebra on the real
coordinate ring $\mathbb R[T^{\mathbb C}\big/W] \cong \mathbb
R[T^{\mathbb C}]^W$ of the adjoint quotient is simply obtained by
taking invariants. This Poisson algebra is somewhat more easily
described in terms of the \emph{complexification\/} $\mathbb
R[T^{\mathbb C}\big/W]_{\mathbb C}$ of the real coordinate ring of
$T^{\mathbb C}\big/W$---this complexification is \emph{not\/} the
complex coordinate ring of $T^{\mathbb C}\big/W$, though. An
explicit description of the induced Poisson structure on the
complexification  of the real coordinate ring of $T^{\mathbb
C}\big/W$ has been worked out in \cite{adjoint}.

For $K=\mathrm{SU}(n)$, the Weyl group $W$ is the symmetric group
$S_n$ on $n$ letters, and the quoted description of the ring
$\mathbb R[T^{\mathbb C}\big/W]_{\mathbb C}$, together with the
Poisson structure,
 takes the following
form: As before, consider $\mathrm{SU}(n)$ as a subgroup of
$\mathrm{U}(n)$ as usual and, accordingly, identify $T^{\mathbb
C}$ with the subgroup of  $(\mathbb C^*)^n$ (taken as a maximal
torus of $\mathrm{U}(n)^{\mathbb C}=\mathrm{GL}(n,\mathbb C)$)
given by $z_1 \cdot\ldots\cdot z_n= 1$ where $z_1$, \dots, $z_n$
are the obvious coordinates on $(\mathbb C^*)^n$. For $r \geq 0$
and $s \geq 0$ such that $1 \leq r + s \leq n$, let
\begin{equation}
\sigma_{(r,s)}(z_1,\dots, z_n,\overline z_1, \dots, \overline z_n)
\label{elementary}
\end{equation}
be the $S_n$-orbit sum of the monomial $ z_1\cdot\ldots\cdot z_r
\overline z_{r+1}\cdot\ldots\cdot \overline z_{r+s}$; the function
$\sigma_{(r,s)}$ is manifestly $S_n$-invariant, and we refer to a
function of this kind as an \emph{elementary bisymmetric
function\/}. In a degree $m$, $1 \leq m \leq n$, the construction
yields the $m+1$ bisymmetric functions $\sigma_{(m,0)}$,\
$\sigma_{(m-1,1)}$, \dots, $\sigma_{(0,m)}$ whence in degrees at
most equal to $n$  it yields altogether $\frac {n(n+3)}2$
elementary bisymmetric functions. According to a classical result,
the ring $\mathbb Q[ z_1,\dots, z_n,\overline z_1, \dots,
\overline z_n]^{S_n}$ of bisymmetric functions (over the
rationals) is generated by the elementary bisymmetric functions.
This implies that, as a ring, $\mathbb R[T^{\mathbb
C}\big/S_n]_{\mathbb C}$ is generated by the $\frac {n(n+3)}2-2$
elementary bisymmetric functions $\sigma_{(r,s)}$ for $0 \leq r
\leq n$ and $0 \leq s \leq n$ such that $1 \leq r + s\leq n$ and
$(r,s) \ne (n,0)$ and $(r,s) \ne (0,n)$, subject to a certain
system of $\frac {n(n-1)}2$ relations. See \cite{adjoint} for
details. For intelligibility we note that the conditions $(r,s)
\ne (n,0)$ and $(r,s) \ne (0,n)$ mean that, in accordance with the
discussion in Section \ref{SL} above, in the ring $\mathbb
R[T^{\mathbb C}\big/S_n]_{\mathbb C}$,
 $\sigma_{(n,0)} =1$ and
$\sigma_{(0,n)} =1$.

The stratified symplectic Poisson structure is more easily
described in terms of  another system of multiplicative generators
for the complexification of the real coordinate ring of the
adjoint quotient; this system is equivalent to the above one and
is obtained in the following fashion: For $r \geq 0$ and $s \geq
0$ such that $1 \leq r + s$, let
\begin{equation}
\tau_{(r,s)}(z_1,\dots, z_n,\overline z_1, \dots, \overline z_n) =
\sum_{j=1}^n z^r_j\overline z^s_j;
\end{equation}
such a function $\tau_{(r,s)}$ is manifestly $S_n$-invariant, and
we refer to a function of this kind as a \emph{bisymmetric power
sum function\/}. As a ring, $\mathbb R[T^{\mathbb
C}\big/S_n]_{\mathbb C}$ is generated by the $\frac {n(n+3)}2$
bisymmetric power sum functions $\tau_{(r,s)}$ for $r \geq 0$ and
$s \geq 0$ such that $1 \leq r + s\leq n$ as well, subject to
$\frac {n(n-1)}2+2$ relations; see \cite{adjoint} for details.
Suffice it to mention at this stage that rewriting the relation
$\sigma_{(n,0)} =1$ in terms of the power sums $\tau_{(j,0)}$ of
the variables $z_1$, \dots, $z_n$ and the relation $\sigma_{(0,n)}
=1$ in terms of the power sums $\tau_{(0,j)}$  of the variables
$\overline z_1$, \dots, $\overline z_n$ where $1 \leq j \leq n$,
we can express the generator $\tau_{(n,0)}$ as a polynomial in the
$\tau_{(j,0)}$'s with $1 \leq j < n$ and, likewise, we can express
the generator $\tau_{(0,n)}$ as a polynomial in the
$\tau_{(0,j)}$'s with $1 \leq j < n$. This procedure reduces the
above system to one with $\frac {n(n+3)}2-2$ generators, subject
to $\frac {n(n-1)}2$ relations.

By a result in \cite{adjoint}, the Poisson brackets among the
multiplicative generators $\tau_{(j,k)}$ {\rm ($0\leq j \leq n, 0
\leq k \leq n, 1 \leq j+k \leq n$)}  are given by the formulas
\begin{equation}
\frac i2\{\tau_{(j_1,k_1)},\tau_{(j_2,k_2)}\} =
(j_1k_2-j_2k_1)\tau_{(j_1+j_2,k_1+k_2)}. \label{poisson}
\end{equation}
When $j_1+j_2+k_1+k_2 >n$, the right-hand side of \eqref{poisson}
is here to be rewritten as a polynomial in terms of the
multiplicative generators $\tau_{(j,k)}$, for $0\leq j \leq n, 0
\leq k \leq n, 1 \leq j+k \leq n$. Thus we see that, as a
\emph{stratified K\"ahler space\/}, the adjoint quotient of
$\mathrm{SL}(n,\mathbb C)$ is considerably more complicated than
just as a complex algebraic variety; indeed, as a complex
algebraic variety, this quotient is just a copy of
$(n-1)$-dimensional affine complex space.

We will now spell out explicitly the stratified K\"ahler structure
for $n=2$ and $n=3$. Consider first the case where $n=2$, so that
$K=\mathrm{SU}(2)$ and $K^{\mathbb C}=\mathrm{SL}(2,\mathbb C)$.
The standard maximal torus $T\cong S^1$ of $\mathrm{SU}(2)$
consists of the diagonal matrices
$\mathrm{diag}(\zeta,\zeta^{-1})$ where $|\zeta|=1$, and the
standard maximal torus $T^{\mathbb C}\cong \mathbb C^*$ of
$\mathrm{SL}(2,\mathbb C)$ consists of the diagonal matrices
$\mathrm{diag}(\zeta,\zeta^{-1})$ where $\zeta\ne 0$.
 In
view of the discussion in Section \ref{SL}, complex algebraically,
the categorical quotient $K^{\mathbb C}\big/\big/K^{\mathbb C}$
amounts to the space $T^{\mathbb C}\big /S_2 \cong \mathbb C$ of
orbits relative to the action of the Weyl group $S_2$ on $\mathbb
C^* \cong T^{\mathbb C}$, and this orbit space is realized as the
target of the holomorphic map
\begin{equation}
\chi \colon \mathbb C^* \longrightarrow \mathbb C,\ \chi(z) = z+
z^{-1}. \label{chi}
\end{equation}
Thus $Z = z+ z^{-1}$ may be taken as a holomorphic coordinate on
the adjoint quotient of $\mathrm{SL}(2,\mathbb C)$.

In view of the above remarks, the complexification $\mathbb
R[T^{\mathbb C}\big/S_2]_{\mathbb C}$ of the \emph{real\/}
coordinate ring $\mathbb R[T^{\mathbb C}\big/S_2]$ of the adjoint
quotient $T^{\mathbb C}\big/S_2$ under discussion is generated by
the three elementary bisymmetric functions
\[
\sigma_1 =\sigma_{(1,0)}= z+z^{-1},\ \overline \sigma_1 =
\sigma_{(0,1)}=\overline z+ \overline z^{-1}, \
 \sigma =\sigma_{(1,1)} =\frac z{\overline z}+\frac {\overline z}{z},
\]
subject to the single defining relation
\[
(\sigma_1^2-4)(\overline \sigma_1^2-4) = (\sigma_1 \overline
\sigma_1 -2\sigma)^2.
\]
See  \cite{adjoint} for details. Hence the real coordinate ring
$\mathbb R[T^{\mathbb C}\big/S_2]$ of the adjoint quotient
$T^{\mathbb C}\big/S_2$ under discussion is generated by the three
functions
\[
X = x + \frac{x}{r^2}, \ Y = y - \frac{y}{r^2},\ \tau =
\frac{2-\sigma}4 = \frac{y^2}{r^2},
\]
where $z=x+iy$, $Z=X+iY$, and $x^2 + y ^2 = r^2$, subject to the
relation
\begin{equation}
Y^2 = (X^2 + Y^2 + 4 (\tau -1)) \tau. \label{relation}
\end{equation}
The obvious inequality $\tau \geq 0$ brings the semialgebraic
nature of the adjoint quotient to the fore. More details
concerning the semialgebraic structure may be found in
\cite{adjoint}. Moreover, the Poisson bracket
$\{\,\cdot\,,\,\cdot\,\}$ on $\mathbb R[T^{\mathbb C}\big/S_2]$ is
given by
\[
\{X,Y\} = X^2 + Y ^2 + 4(2 \tau -1)
\]
together with
\[
\{X,\tau\} = 2(1- \tau)Y,\ \{Y,\tau\} = 2 \tau X.
\]
The resulting complex algebraic stratified K\"ahler structure is
\emph{singular\/} at the two points  $-2$ and $2$. These are
solutions of the discriminant equation \eqref{equation13}, and the
Poisson structure vanishes at these points; furthermore, at these
two points, the function $\tau$ is \emph{not\/} an ordinary smooth
function of the variables $X$ and $Y$. Indeed, solving
\eqref{relation} for $\tau$, we obtain
$$
\tau = \frac 12 \sqrt{Y^2+ \frac {(X^2 + Y^2-4)^2}{16}} -
\frac{X^2 + Y^2 -4}{8} ,
$$
whence, at $(X,Y)=(\pm 2,0)$, $\tau$ is not smooth as a function
of the variables $X$ and $Y$. Away from these two points, the
Poisson structure is symplectic. We refer to the adjoint quotient
under discussion as the \emph{exotic plane with two vertices\/}.

This exotic plane with two vertices admits the following function
theoretic interpretation: The composite
\[
\chi \circ\mathrm{exp} \colon \mathfrak t^{\mathbb C} \cong
\mathbb C \longrightarrow \mathbb C \cong T^{\mathbb C}\big/S_2
\]
equals the holomorphic function $2\mathrm{cosh}$. Consequently, under the
composite
\[
S^1 \times \mathbb R \longrightarrow \mathbb C^* \longrightarrow
\mathbb C
\]
of the polar map with $\chi$, the family of circles $S_t^1
=\{(\eta,t); \eta \in S^1\}$ ($t \in \mathbb R$) in $S^1 \times
\mathbb R$ (each circle of the family being parallel to the zero
section) goes to the family of curves
\[
\eta \mapsto \mathrm e^t \eta +\mathrm e^{-t}\eta^{-1},\quad \eta
\in S^1,
\]
which are ellipses for $t\ne 0$; and the family of lines
$L_{\eta}=\{(\eta,t);t \in \mathbb R\}$ ($\eta \in S^1$) in the
tangent directions goes to the family of curves
\[
t \mapsto \mathrm e^t \eta +\mathrm e^{-t}\eta^{-1},\quad t \in
\mathbb R,
\]
which are hyperbolas for $\eta \ne \pm 1$. The image of the circle
$S_0^1$ is a double line segment between $-2$ and $2$, the image
of the line $L_{1}$ is a double ray, emanating from $2$ into the
positive real direction, and the image of the line $L_{-1}$ is a
double ray, emanating from $-2$ into the negative real direction.
The two families are orthogonal and have the two singular points
$-2$ and $2$ as focal points. This is of course well known and
entirely classical. In the cotangent bundle picture, the double
line segment between $-2$ and $2$ is the adjoint quotient of the
base $K=\mathrm {SU}(2)$, indeed this line segment is exactly the
fundamental alcove  of $\mathrm {SU}(2)$, the family of hyperbolas
which meet this line segment between $-2$ and $2$ constitutes a
cotangent bundle on this orbit space with the two singular points
removed, and the plane, i.~e. adjoint quotient of $K^{\mathbb C}$,
with the two rays emanating from $-2$ and $2$ removed, is the
total space of this cotangent bundle; furthermore, the cotangent
bundle symplectic structure is precisely that which corresponds to
the reduced Poisson structure. However the cotangent bundle
structure does not extend to the entire adjoint quotient of
$\mathrm T^*K\cong K^{\mathbb C}$. In particular, this
interpretation visualizes the familiar fact that, unless there is
a single stratum, the strata arising from cotangent bundle
reduction are {\it not\/} cotangent bundles on strata of the orbit
space of the base space.

Consider now the case where $K=\mathrm {SU}(3)$, with maximal
torus $T$ diffeomorphic to the product $S^1 \times S^1$ of two
copies of the circle group. Complex algebraically, the map
\eqref{steinberg2} for $n=3$ comes down to the map
\begin{equation}
(-\sigma_1,\sigma_2) \colon T^{\mathbb C} \longrightarrow \mathbb
C^2, \label{steinberg3}
\end{equation}
and this map induces a
complex algebraic isomorphism
from the adjoint quotient
$T^{\mathbb C}\big/S_3$ of $\mathrm {SL}(3,\mathbb C)$ onto a copy
$\mathbb C^2$ of 2-dimensional complex affine space; as in Section
\ref{SL} above, we will take $a_1=-\sigma_1$ and $a_2=\sigma_2$ as
complex coordinates on the adjoint quotient. The complement of the
top stratum $\mathcal D^{(o,\mathrm{SU}(3))}_3$ is now the complex
affine curve $\mathcal D^{\mathrm{SU}(3)}_2= \mathcal
D^{\mathrm{SU}(3)}_{(2,1)}$ in $\mathbb A^{3,1}_{\mathrm{coef}}$
given by the equation \eqref{equation14}. This curve is plainly
parametrized by the restriction
\begin{equation}
\mathbb C^* \longrightarrow \mathbb C \times \mathbb C, \quad z
\mapsto (2z + z^{-2}, z^2 + 2z^{-1}) \label{steinberg4}
\end{equation}
of \eqref{steinberg3} to the diagonal, and this holomorphic curve
parametrizes the {\it closure\/} $\mathcal D^{\mathrm{SU}(3)}_2$
of the stratum $\mathcal D^{(o,\mathrm{SU}(3))}_2$. This curve has
the three (complex) singularities $(3,3)$, $3(\eta,\eta^2)$,
$3(\eta^2,\eta)$, where $\eta^3 = 1, \eta \ne 1$. These points are
the images of the (conjugacy classes) of the three central
elements under \eqref{steinberg3}; as complex curve singularities,
these singularities are cuspidal. These three points constitute
the stratum $\mathcal D^{\mathrm{SU}(3)}_3$.

The real hypersurface written above, for general $n$, as $\mathcal
H_n$, cf. \eqref{hypersurface}, now comes down to the curve $
\mathcal H_2$ in $\mathbb R^2$; here $\mathbb R^2$ is embedded
into $\mathbb C^2$ as the real affine space of real points of
$\mathbb C^2$ in the obvious fashion. Since, for a complex number
$z$ with $|z|=1$, $2\overline z + \overline z^{-2}= z^2 +
2z^{-1}$, the restriction of \eqref{steinberg4} to the {\it
real\/} torus $T\subseteq T^{\mathbb C}$ yields the parametrized
real curve
\begin{equation}
S^1 \longrightarrow \mathbb R^2,\quad \mathrm e^{i\alpha} \mapsto
(u(\alpha),v(\alpha))\in \mathbb R^2, \label{steinberg5}
\end{equation}
where $u(\alpha) + i v(\alpha) = 2 \mathrm e^{i\alpha} + \mathrm
e^{-2i\alpha}$, and this parametrizes precisely the curve
$\mathcal H_2$. This curve  is a \emph{hypocycloid\/}, as noted in
\cite{chakirus} (Section 5), and the real orbit space $\mathrm
{SU}(3)\big/\mathrm {SU}(3)$ relative to conjugation, realized
within the adjoint quotient $T^{\mathbb C}\big/S_3\cong \mathcal
D^{\mathrm{SU}(3)}_1$ of $\mathrm {SL}(3,\mathbb C)$, amounts to
the compact region in $\mathbb R^2$ enclosed by this hypocycloid.

As a complex algebraic stratified K\"ahler space, the adjoint
quotient looks considerably more complicated. Indeed, cf.
\cite{adjoint}, the complexification $\mathbb R[T^{\mathbb
C}\big/S_3]_{\mathbb C}$ of the real coordinate ring $\mathbb
R[T^{\mathbb C}\big/S_3]\cong \mathbb R[T^{\mathbb C}]^{S_3}$ of
the adjoint quotient $T^{\mathbb C}\big/S_3$, viewed as a real
semialgebraic space, is generated by the seven functions $\sigma_1
=\sigma_{(1,0)}$, $\overline \sigma_1=\sigma_{(0,1)}$,  $\sigma_2
=\sigma_{(2,0)}$, $\overline \sigma_2=\sigma_{(0,2)}$,
$\sigma=\sigma_{(1,1)}$,  $\rho =\sigma_{(2,1)}$, $\overline
\rho=\sigma_{(1,2)}$, subject to the following three relations:
\begin{align*}
(\sigma_1^2-4 \sigma_2)(\overline \sigma_1^2-4 \overline \sigma_2)
&= (\sigma_1 \overline \sigma_1 -2\sigma)^2 + 2\rho \overline
\sigma_1 + 2\overline\rho \sigma_1
\\
D_3(1,-\sigma_1,\sigma_2, - 1) \overline \sigma_2 &=(9 +
\sigma_2^3 - 4 \sigma_1 \sigma_2) \overline \sigma_1^2
\\
{}& \quad +(4 \sigma_1^2  -3 \sigma_2 - \sigma_1 \sigma_2^2)
\overline \sigma_1 \sigma
\\
{}& \quad +(6 \sigma_1  - \sigma_2^2) \overline \sigma_1 \rho
+(\sigma_2^2-3\sigma_1  )  \sigma^2
\\
{}& \quad +(9  - \sigma_1 \sigma_2) \sigma \rho +(\sigma_1^2 - 3
\sigma_2) \rho^2
\\
D_3(1,-\sigma_1,\sigma_2, - 1) &=
 \overline \sigma_1^3
-\sigma_2 \overline \sigma_1^2 \sigma + (\sigma_2^2 - 2 \sigma_1)
\overline \sigma_1^2 \rho + \sigma_1 \overline \sigma_1 \sigma^2
\\
{}&\quad -((\sigma_1^2 - 2 \sigma_2)\sigma_1 - \sigma_1^3 + 3
\sigma_1 \sigma_2 - 3) \overline \sigma_1 \sigma \rho
\\
{}&\quad -\sigma^3 +(\sigma_1^2 - 2 \sigma_2) \overline \sigma_1
\rho^2
\\
{}&\quad+ \sigma_2 \sigma^2 \rho - \sigma_1 \sigma \rho^2 +
\rho^3.
\end{align*}
We note that an explicit expression for $D_3(1,-\sigma_1,\sigma_2,
- 1)$ is given by \eqref{cubic}. The formula \eqref{poisson}
 yields the stratified symplectic Poisson
structure on the complexification of the real coordinate ring of
$T^{\mathbb C}\big/S_3$ in terms of the nine  generators
$\tau_{(1,0)}$, $\tau_{(2,0)}$, $\tau_{(3,0)}$, $\tau_{(0,1)}$,
$\tau_{(0,2)}$, $\tau_{(0,3)}$, $\tau_{(1,1)}$, $\tau_{(1,2)}$,
$\tau_{(2,1)}$. This Poisson structure has rank 4 on the top
stratum $\mathcal D^{(o,\mathrm{SU}(3))}_1$, rank 2 on the stratum
$\mathcal D^{(o,\mathrm{SU}(3))}_2$, and rank zero at the three
points of the stratum $\mathcal D^{\mathrm{SU}(3)}_1$, that is, at
the three cusps of the complex affine curve $\mathcal
D^{\mathrm{SU}(3)}_2= \mathcal D^{\mathrm{SU}(3)}_{(2,1)}$ in
$\mathbb A^{3,1}_{\mathrm{coef}}$ given by the equation
\eqref{equation14}. The requisite inequalities which encapsulate
the semialgebraic structure may be found in \cite{adjoint}.

For $n \geq 4$, the map \eqref{steinberg2} identifies the adjoint
quotient $\mathcal D^{\mathrm{SU}(n)}_1$ of $\mathrm
{SL}(n,\mathbb C)$ with complex affine $(n-1)$-space $\mathbb
C^{n-1}$, and the resulting stratified K\"ahler structure on the
adjoint quotient can be described in a way similar to that for the
low dimensional cases where $n=2$ and $n=3$; the details get more
and more involved, though. A function theoretic interpretation
extending the interpretation spelled out above for the case where
$n=2$ is available as well. Indeed, let $n\geq 2$, let $T\cong
(S^1)^{(n-1)}$ be the standard maximal torus in $\mathrm {SU}(n)$,
let $T^{\mathbb C}\cong (\mathbb C^*)^{n-1}$ be its
complexification, and consider  the composite
\begin{equation}
T \times \mathfrak t \cong T \times \mathbb R^{n-1}
\longrightarrow T^{\mathbb C}\longrightarrow \mathbb C^{n-1}
\label{cosh}
\end{equation}
of the polar map with $\chi$. Now the images in the top stratum of
the adjoint quotient $\mathcal D^{\mathrm{SU}(n)}_1\cong \mathbb
C^{n-1}$ of the leaves of the horizontal foliation of the total
space $T \times \mathbb R^{n-1} \cong \mathrm T T$ of the tangent
bundle of $T$ are smooth submanifolds and generalize the family of
ellipses for the case where $n=2$ and, likewise, the images in the
top stratum of the adjoint quotient $\mathbb C^{n-1}$ of the
leaves of the vertical foliation of the total space $T \times
\mathbb R^{n-1} \cong \mathrm T T$ of the tangent bundle of $T$
are smooth manifolds and generalize the family of hyperbolas for
the case where $n=2$; however, images of such leaves which meet a
lower stratum involve singularities, and the geometry of the
situation can be described in terms of complex focal points,
complex focal lines, etc.: The focal points are the $n$ points of
the bottom stratum $\mathcal D^{\mathrm{SU}(n)}_n$, a focal curve
is given by the next stratum $\mathcal D^{\mathrm{SU}(n)}_{n-1}$,
etc. The composite of \eqref{cosh} with the obvious map from
$\mathbb C^{n-1}$ to $T \times \mathbb R^{n-1}$ generalizes the
holomorphic hyperbolic cosine function for the case where $n=2$,
and the function theoretic interpretation extending the
interpretation for the case where $n=2$ is thus immediate.

The above observations involving the singular cotangent bundle
projection map for the special case where $n=2$ extend as follows:
The induced projection from the adjoint quotient $\mathcal
D^{\mathrm{SU}(n)}_1\cong\mathbb C^{n-1}$ to the real orbit space
$\mathrm {SU}(n)\big/ \mathrm {SU}(n)$ relative to conjugation is
a singular cotangent bundle projection map in an obvious manner.
Indeed, the map \eqref{steinberg2} induces a map from the
\emph{real} torus $T$ to the real affine space $\mathbb R^{n-1}$
of real points of the complex adjoint quotient $\mathrm C^{n-1}$;
just as for the cases where $n=2$ or $n=3$,  the restriction of
this map to the subspace $T_1$ of the real torus $T$ having
non-trivial stabilizer, i.~e. points that are non-regular as
points of $K=\mathrm {SU}(n)$, has as its image a real closed
hypersurface in $\mathbb R^{n-1}$, and the real orbit space
$\mathrm {SU}(n)\big/ \mathrm {SU}(n)$ relative to conjugation is
realized in $\mathbb R^{n-1}$ as the semialgebraic space enclosed
by this hypersurface. The resulting semialgebraic space is
manifestly homeomorphic to an $(n-1)$-simplex or, equivalently, to
the standard $(n-1)$-alcove for $\mathrm {SU}(n)$.

Finally we note that this kind of interpretation is available for
the closure $\mathcal D^{\mathrm{SU}(n)}_{\nu}$ of each stratum of
the adjoint quotient, each such space being a complex algebraic
stratified K\"ahler space; moreover, for $\nu' \preceq \nu$, the
injection map from $\mathcal D^{\mathrm{SU}(n)}_{\nu'}$ into
$\mathcal D^{\mathrm{SU}(n)}_{\nu}$ is compatible with all the
structure. For example, when $n=3$, $\mathcal
D^{\mathrm{SU}(3)}_2$ is the complex affine curve in the adjoint
quotient  $\mathcal D^{\mathrm{SU}(3)}_1 \cong \mathbb C^2$
parametrized by \eqref{steinberg4}. This parametrization plainly
factors through $\mathbb C^*\big/S_2$ and hence induces an
embedding
\[
\mathcal D^{\mathrm{SU}(2)}_1 \longrightarrow \mathcal
D^{\mathrm{SU}(3)}_1
\]
of the adjoint quotient $\mathcal D^{\mathrm{SU}(2)}_1$ of
$\mathrm{SL}(2,\mathbb C)$ into the adjoint quotient $\mathcal
D^{\mathrm{SU}(3)}_1$ of $\mathrm{SL}(3,\mathbb C)$ which
identifies $\mathcal D^{\mathrm{SU}(2)}_1$ with $\mathcal
D^{\mathrm{SU}(3)}_2$ as complex algebraic stratified K\"ahler
spaces. Hence the embedding is one of complex algebraic stratified
K\"ahler spaces, that is, it is compatible with all the structure.
In particular, the \lq\lq focal\rq\rq\ geometries correspond.

\section{Relationship with the spherical pendulum}

In this section we will identify the adjoint quotient of
$\mathrm{SL}(2,\mathbb C)$ with the reduced phase space of a
spherical pendulum constrained to move with angular momentum zero
about the third axis, so that it amounts to a \emph{planar}
pendulum.

The unreduced phase space of the spherical pendulum is the total
space $\mathrm T^*S^2$ of the cotangent bundle of the ordinary
2-sphere $S^2$ in 3-space $\mathbb R^3$ centered at the origin. By
means of the standard inner product, we identify the tangent
bundle of $\mathbb R^3$ with its cotangent bundle;  accordingly,
we identify $\mathrm T^*S^2$ with the total space $\mathrm T S^2$
of the tangent bundle of $S^2$, and we realize the space $\mathrm
T S^2$ within $ \mathbb R^3 \times \mathbb R^3\cong \mathrm T
\mathbb R^3$ in the standard fashion. The map
\begin{equation}
\mathbb C^* \longrightarrow \mathbb R^3 \times \mathbb R^3, \
\mathrm e^t \mathrm e^{i\varphi} \longmapsto(0,\sin \varphi ,\cos
\varphi,0,-t\cos \varphi,t\sin \varphi), \label{cylinder1}
\end{equation}
plainly goes into the subspace $\mathrm T S^2$ of $\mathbb
R^3\times\mathbb R^3$ and induces a homeomorphism from the
$S_2$-orbit space $\mathbb C^*\big/ S_2$ onto the angular momentum
zero reduced phase  space $V_0=\mu^{-1}(0)\big/ S^1$ of the
spherical pendulum. Here the non-trivial element of the cyclic
group $S_2$ with two elements acts on $\mathbb C^*$ by inversion,
and
\begin{equation*} \mu\colon \mathrm T S^2 \longrightarrow
\mathrm{Lie}(S^1)^*\cong \mathbb R
\end{equation*}
refers to the momentum mapping for the $S^1$-action on $\mathrm T
S^2$ given by rotation about the third axis; actually this is just
ordinary angular momentum given by
\[
\mu(q_1,q_2,q_3,p_1,p_2,p_3) =q_1 p_2 - q_2 p_1
\]
where $(q_1,q_2,q_3,p_1,p_2,p_3)$ are the obvious coordinates on
$\mathbb R^3 \times \mathbb R^3$.
 Under the identification of $\mathbb C^*$ with the
complexification $T^{\mathbb C}$ of the maximal torus $T$ of
$\mathrm{SU}(2)$, the quotient $\mathbb C^*\big/ S_2$ gets
identified with  the adjoint quotient $T^{\mathbb C}\big/ S_2$ of
$\mathrm{SL}(2,\mathbb C)$ whence the map \eqref{cylinder1}
induces a homeomorphism from the adjoint quotient $T^{\mathbb
C}\big/ S_2$ of $\mathrm{SL}(2,\mathbb C)$ onto the reduced phase
space $V_0$ for the spherical pendulum. Moreover,  a little
thought reveals that the map
\begin{equation*}
\Theta\colon \mathbb C^* \longrightarrow \mathbb R^3, \
\Theta(\mathrm e^t \mathrm e^{i\varphi})= (\cos \varphi, t \sin
\varphi, t^2)
\end{equation*}
into $\mathbb R^3$ induces a homeomorphism from $\mathbb C^*\big/
S_2$ onto the real semialgebraic subspace $M_0$ of $\mathbb R^3$
which, in terms of the coordinates $(u,v,w)$, is
 given by
\begin{equation*}
w(1-u^2)-v^2 = 0,\ |u| \leq 1,\ w - v^2 \geq 0.
\end{equation*}
The  space $M_0$ has received considerable attention in the
literature: This space is the \lq\lq canoe\rq\rq, cf.
\cite{armcusgo}, \cite{cushbate} (pp.~148 ff.) and
\cite{lermonsj}, the two singular points being absolute equilibria
for the spherical pendulum.  Thus the map $\Theta$ identifies the
adjoint quotient $T^{\mathbb C}\big/ S_2$ of
$\mathrm{SL}(2,\mathbb C)$ and hence the reduced phase space $V_0$
for the spherical pendulum with the \lq\lq canoe\rq\rq\ $M_0$.

We have already seen that the holomorphic map $\chi$ from $\mathbb
C^*$ to $\mathbb C$ given by $\chi(z) = z + z^{-1}$ induces a
complex algebraic isomorphism from
 the adjoint quotient $T^{\mathbb C} \big/
S_2$ onto a copy of the complex line $\mathbb C$. It is
instructive to realize $\chi$ on the homeomorphic image $V_0$ or,
equivalently, on the \lq\lq canoe\rq\rq\  $M_0$. We can even do
better: Let $\alpha$ be the real analytic function of the real
variable $t$ given by the power series $\alpha(t) =
\sum_{j=0}^{\infty} \frac {t^j}{(2j)!}$ and, likewise, let $\beta$
be the real analytic function of the real variable $t$ given by
the power series $\beta(t) = \sum_{j=0}^{\infty} \frac
{t^j}{(2j+1)!}$. We note that, by construction,
\[
\alpha(t^2) = \mathrm{cosh}(t), \ t\beta(t^2) = \mathrm{sinh}(t) .
\]
Introduce the real analytic function
\[
\Phi\colon  \mathbb R^3 \longrightarrow \mathbb C, \ \Phi(u,v,w) =
u\alpha(w) + i v\beta(w) ,(u,v,w) \in \mathbb R^3.
\]
Since the composite map $\Phi \circ \Theta \circ \mathrm{exp}$
from $\mathbb C$ to $\mathbb C$ coincides with $2 \mathrm{cosh}$,
the composite $\Phi \circ \Theta$ coincides with $\chi$. In
particular, the homeomorphism from the \lq\lq canoe\rq\rq\  $M_0$
onto the complex line given by the restriction of $\Phi$
\emph{flattens out\/} the \lq\lq canoe\rq\rq. It is interesting to
note that this flattening out is accomplished in the category of
real semianalytic spaces, not in that of real semialgebraic
spaces, even though the \lq\lq canoe\rq\rq\ and the complex line
are real semialgebraic spaces; presumably the two cannot be
identified in the category of real semialgebraic spaces.

Since the composite $\Phi \circ \Theta\circ\mathrm{exp}$ equals
twice the holomorphic hyperbolic cosine function, the present
discussion provides a geometric interpretation of the adjoint
quotient of $\mathrm{SU}(2)^{\mathbb C}\cong \mathrm{SL}(2,\mathbb
C)$ similar to the function theoretic one given earlier, but in
terms of the \lq\lq canoe\rq\rq. In particular, the two focal
points in the earlier interpretation now correspond to the two
singular points of the \lq\lq canoe\rq\rq\ which, in turn,
correspond to the absolute equilibria of the spherical pendulum.

While the spherical pendulum has been extensively studied, the
present complex analytic interpretation seems to be new.

\section*{Acknowledgments}

This paper was written during a stay at the
Institute for theoretical physics at the
University of Leipzig. This stay was
made possible by the Deutsche Forschungsgemeinschaft in the framework
of a Mercator-professorship, and
I wish to express my gratitude to the
Deutsche Forschungsgemeinschaft.
It is a pleasure to acknowledge the stimulus of conversation with
G. Rudolph and M. Schmidt at Leipzig.

\bigskip

\noindent Universit\'e des Sciences et Technologies de Lille, UFR
de Math\'ematiques, CNRS-UMR 8524
\\ 59655 VILLENEUVE D'ASCQ, C\'edex, France
\\ and
\\ Institute
for Theoretical Physics, Universit\"at Leipzig\\
04109 LEIPZIG, Germany \\
{Johannes.Huebschmann@math.univ-lille1.fr}
\end{document}